\documentclass[12pt]{article}

\usepackage{times}

\def\x#1{\def\@@eqncr{\let\@tempa\relax
  \ifcase\@eqcnt \def\@tempa{& & &}\or \def\@tempa{& &}
   \else \def\@tempa{&}\fi & #1\cr}}

\def\ra{\rightarrow}

\def\ss{\subseteq}

\def\Re{\hbox{\rm Re}\,}
\def\Im{\hbox{\rm Im}\,}

 \def\HollowBox #1#2{{\dimen0=#1 \advance\dimen0 by -#2
       \dimen1=#1 \advance\dimen1 by #2
        \vrule height #1 depth #2 width #2
        \vrule height 0pt depth #2 width #1
        \llap{\vrule height #1 depth -\dimen0 width \dimen1}%
       \hskip -#2
       \vrule height #1 depth #2 width #2}}
 \def\BoxOpTwo{\mathord{\HollowBox{6pt}{.4pt}}\;}

\def\endpf{\hfill $\BoxOpTwo$}

\font\teneufm=eufm10
\font\seveneufm=eufm7
\font\fiveeufm=eufm5
\newfam\eufmfam
\textfont\eufmfam=\teneufm
\scriptfont\eufmfam=\seveneufm
\scriptscriptfont\eufmfam=\fiveeufm

\newfam\msbfam
\font\tenmsb=msbm10  scaled \magstep1 \textfont\msbfam=\tenmsb
\font\sevenmsb=msbm7  scaled \magstep1 \scriptfont\msbfam=\sevenmsb
\font\fivemsb=msbm5  scaled \magstep1  \scriptscriptfont\msbfam=\fivemsb
\def\Bbb{\fam\msbfam \tenmsb}

\def\RR{{\Bbb R}}
\def\CC{{\Bbb C}}

\newtheorem{theorem}{Theorem}
\newtheorem{corollary}[theorem]{Corollary}
\newtheorem{proposition}[theorem]{Proposition}
\newtheorem{lemma}[theorem]{Lemma}

\def\rnp{\RR^{N+1}_+}
\def\dbar{\overline{\partial}}

\def\Re{\hbox{Re}\,}
\def\Im{\hbox{Im}\,}

\begin{document}

\null \vskip-.5in \hfill {\sc Preliminary Draft} \\
\rightline{\it September 11, 2005}
\vspace*{.3in}

\begin{center}
\Large \bf The Boundary Behavior of Holomorphic Functions:
\smallskip \\
\Large Global and Local Results
\medskip \\
{\normalsize \rm by Steven G. Krantz}
\end{center}
\vskip.25in

\begin{quote}
\noindent {\footnotesize \bf Abstract:} \footnotesize \sl We develop
a new technique
for studying the boundary limiting behavior of a holomorphic
function on a domain $\Omega$---both in one and several
complex variables. The approach involves two new localized
maximal functions.

\ \ \ \ As a result of this methodology, theorems of Calder\'{o}n type
about local boundary behavior on a set of positive measure may
be proved in a new and more natural way.

\ \ \ \ We also study the question of nontangential boundedness (on a
set of positive measure) versus admissible boundedness. Under
suitable hypotheses, these two conditions are shown to be equivalent.
\end{quote}

\setcounter{section}{-1}

\section{Introduction}

The first theorem about the boundary limiting behavior of holomorphic
functions was proved by P. Fatou in his thesis in 1906 [FAT]. He used Fourier
series techniques to show that if $f$ is a bounded, holomorphic
function on the unit disc $D \ss \CC$ (i.e., $f \in H^\infty(D)$), then
$$
\lim_{r \ra 1^-} f(r e^{i\theta})
$$
exists for almost every $\theta \in [0,2\pi)$. Furthermore,
for $\alpha > 1$ and $P \in \partial D$, we define
$$
\Gamma_\alpha(P) = \{\zeta \in D: |z - P| < \alpha (1 - |z|)\} \, .
$$
This is the {\it nontangential} or {\it Stolz} approach region.
Then Fatou showed that, for $f \in H^\infty(D)$ and $\alpha > 1$ fixed,
the limit
$$
\lim_{\Gamma_\alpha(P) \ni z \ra P} f(z)
$$
exists for almost every $P \in \partial D$.

Later on, Privalov [PRI1], [PRI2], Plessner [PLE], and others refined Fatou's result (see [DIK] for
a detailed account of the history). The standard
theorem today is that, if $0 < p \leq \infty$, if $\alpha > 1$ is fixed, and if
$f \in H^p(D)$ (the Hardy space), then
$$
\lim_{\Gamma_\alpha(P) \ni z \ra P} f(z)
$$
exists for almost every $P \in \partial D$. With suitable estimates for the Poisson
kernel (see [KRA8]), one can prove a similar result
(for nontangential convergence) on any bounded domain in $\CC$ with $C^2$
boundary.  The result may be refined further so that it is valid for functions
in the Nevanlinna class (see [GAR] as well as our
Section 9).

A theorem of Littlewood (see [LIT]) shows
that, in a very precise sense, the nontangential approach regions $\Gamma_\alpha$
are the {\it broadest} approach regions through which a theorem of this kind
can be obtained---even for bounded holomorphic functions.

On a bounded domain in $\CC^n$ with $C^2$ boundary, it is classical that a
holomorphic function satisfying suitable growth conditions---say
membership in the Hardy class $H^p$ or even the Nevanlinna class (see Section 9)---will
have nontangential boundary
limits almost everywhere with respect to $(2n-1)$-dimensional area measure $d\sigma$
on the boundary.  In the present context, ``nontangential'' means
through an approach region of {\it conical shape}:
$$
\Gamma_\alpha(P) = \{z \in \Omega: |z - P| < \alpha \cdot \delta_{\partial \Omega}(P)\} \, ,
$$
with $\delta_{\partial\Omega}(P)$ denoting the Euclidean distance of $P$
to $\partial \Omega$.
Also, if $\Omega = \{z \in \CC^n: \rho(z) < 0\}$, so that $\rho$ is a defining
function for $\Omega$ and if $0 < p < \infty$, then
$$
H^p(\Omega) = \{f \ \hbox{on} \ \Omega: \sup_{0 < \epsilon < \epsilon_0} \int_{\partial \Omega_\epsilon}
|f(\zeta)|^p \, d\sigma(\zeta)^{1/p} \equiv \|f\|_{H^p(\Omega)} < \infty\} \, .
$$
Also $H^\infty(\Omega)$ consists of the bounded holomorphic functions
with the supremum norm.
Here $\epsilon_0$ is a small, positive number and
$\Omega_\epsilon \equiv \{z \in \Omega: \rho(z) = - \epsilon\}$.
Details connected with this definition may be found in [KRA1].
See [TSU], [ZYG1], [ZYG2], [ZYG3]
for historical background of these ideas.
The proof of such a theorem depends once again
on having the appropriate
estimates for the Poisson kernel (see [KRA8], [KRA1]).

It came as quite a surprise when, in 1970, Adam Koranyi [KOR1], [KOR2]
showed that a broader method of approach than
nontangential is valid when the domain in
question is the unit ball in $\CC^n$. To wit, let $B = \{z \in \CC^n:
|z|^2 < 1\}$. For $\alpha > 1$ and $P \in \partial B$, define
$$
{\cal A}_\alpha(P) = \{z \in B: |1 - z \cdot \overline{P}|
  < \alpha (1 - |z|)\} \, .
$$
Here, as is standard, $z \cdot \overline{P} \equiv \sum_j z_j \overline{P}_j$. One
may calculate (see [KRA1]) that the approach region ${\cal A}_\alpha$
has nontangential shape in complex normal directions but parabolic shape in
complex tangential directions. Koranyi's result
is that, if $f \in H^p(B)$, then
$$
\lim_{{\cal A}_\alpha(P) \ni z \ra P} f(z)
$$
exists for $\sigma$-almost every $P \in \partial B$. Koranyi's proof
depends decisively on an analysis of the shape of the
singularity of the Poisson-Szeg\H{o} kernel
$$
{\cal P}(z, \zeta) = \frac{(n - 1)!}{2\pi^n} \frac{(1 - |z|^2)^n}{|1 - z\cdot \overline{\zeta}|^{2n}} .
$$
for the ball (which shape is decidedly different from the shape
of the singularity for the classical Poisson kernel---see
[KRA1], [KRA8]). Put in other words, where the classical results
depend on estimates for the standard Poisson kernel (in
particular, an analysis of its singularity), the new results
of Koranyi required estimates on the Poisson-Szeg\"{o} kernel.

In 1972, E. M. Stein [STE1] showed how to prove a result like
Koranyi's on {\it any} domain in $\CC^n$ with $C^2$ boundary.
His analysis (later refined by Barker---see [BAR] and the discussion in
[KRA1]) avoids the use of canonical kernels, but instead
depends on an analysis of the Levi geometry of the
domain. See also [LEM] for a quite different and
original approach to these matters.

We now realize that Stein's result was an important first step,
but it is far from the optimal result for most domains. More
precisely, the parabolic approach in complex tangential
directions is really only suitable at {\it strongly
pseudoconvex} boundary points. For points of {\it finite type}
$m$, an approach region which has aperture $\Im z_1 = |z'|^m$
(where $z_1$ is the complex normal direction and $z'$ the
remaining complex tangential directions---see the discussion
below) is the right idea.\footnote{It is instructive to
examine the boundary behavior of a holomorphic function a
point $P \in \partial \Omega$ which is {\it strongly
pseudoconcave}. By the {\it Kontinuit\"{a}tssatz} of
multivariable complex function theory (or the Hartogs
extension phenomenon), {\it any} holomorphic function on
$\Omega$ will continue analytically to an entire neighborhood
of $P \in \partial \Omega$. So the correct approach region at
such a boundary point will be unrestricted. This observation
is quite different from the result of Koranyi/Stein.}

But in fact the analysis is much more subtle than that, for it
is not just the {\it type} of the boundary point but the {\it
magnitude of that type} that must play a role. The type and
the magnitude of the type depend semicontinuously on $P \in
\partial \Omega$ when $\Omega \ss \CC^2$. The dependence is
more subtle for $\Omega \ss \CC^n$, $n > 2$. The calculations
in [NSW1], [NSW2] (see also [KRA9]) begin to show how to tame
these rather complicated ideas.

Now let us examine a slightly different direction of these studies.
It is an old result of A. P. Calder\'{o}n [CAL] that if a function
$u$, harmonic on the upper half-
space $\RR^{N+1}_+$, is
nontangentially bounded on a set $E \ss \RR^N \equiv \partial
\RR^{N+1}_+$ of positive $N$-dimensional measure then $u$ has
nontangential boundary limits at almost every point of $E$. In
the important book [STE1], E. M. Stein proved an analogous
result for a holomorphic function on a strongly
pseudoconvex domain $\Omega \ss
\CC^n$ and for admissible boundedness (see [KRA1] as well as
the forthcoming [DIK] for a discussion of both nontangential
and admissible approach regions). Of course it is true
that a holomorphic function on $\Omega \ss \CC^n$
that is nontangentially bounded on a set $E \ss \partial \Omega$ of
positive $(2n-1)$-dimensional measure will have nontangential limits
almost everywhere on $E$---see [CAL], [STE1] and references therein.
So certainly a function that is {\it admissibly bounded} on
a set $E \ss \partial \Omega$ will have nontangential limits
almost everywhere in $E$, just because admissible boundedness is
a stronger condition than nontangential boundedness.  One of the
main points of the present paper is to give a new proof of a fairly
general version of the Calder\'{o}n theorem for admissible
approach regions---one in which the admissible regions have
a geometry that is adapted to the particular domain under study (see Section 8).
Thus the result presented here is more general than that in
[STE1] or [BAR]---just because the approach regions now
fit the Levi geometry.

One of the main points of the present work is to make a comparison
between nontangential behavior of holomorphic functions of several
variables and admissible behavior.  We prove the somewhat surprising
result that if a holomorphic function on the ball $B \subseteq \CC^n$ is
{\it nontangentially bounded} almost everywhere on a set $E \subseteq \partial B$ of positive
measure then it is in fact {\it admissibly} bounded almost everywhere
on $E$.  Discussion, context, and proof appear below.

\section{Nontangential Boundary Behavior Versus \hfill \break
\null \indent Admissible Boundary Behavior}

There are a variety of results in the subject that link, or at least
compare and contrast, the isotropic behavior suggested by
nontangential approach regions with the nonisotropic behavior suggested
by admissible approach regions. Perhaps the first result of this
kind was announced by Stein in [STE2]. The details of
the argument appear in [KRA1]. The result states that
a holomorphic function on a smoothly bounded domain $\Omega$ in $\CC^n$ which
is in a classical Lipschitz function space is in
fact in a stronger nonisotropic function space.
Roughly speaking, such a function is automatically twice
as smooth in tangential directions. This result is in fact
valid for all $0 < \alpha < \infty$.  We now provide some
of the concepts and details pertaining to Stein's result.

For $0 < \alpha < 1$ and $\Omega \ss \RR^N$, we define
$$
\Lambda_\alpha(\Omega) = \{f: \Omega \ra \CC :
 |f(x+h) - f(x)| \leq C |h|^\alpha \ \hbox{for all}\ x, x+h\in \Omega\} \, .
$$
We equip $\Lambda_\alpha$ with the norm
$$
\|f\|_{\Lambda_\alpha} = \sup_{x, x+h \in \Omega \atop h \ne 0}
  \frac{|f(x+h) - f(x)|}{|h|^\alpha} + \|f\|_{L^\infty(\Omega)} \, .
$$

If $\alpha = 1$, then we use the slightly more subtle definition
\begin{eqnarray*}
\Lambda_1(\Omega) & = & \left \{ f: \Omega \ra \CC :
 |f(x + h) + f(x - h) - 2f(x)| \right. \\
\null \qquad \qquad && \hbox{ \ \ \ \ }  \leq \left. C |h| \ \hbox{for all}\ x, x+h, x-h
\in \Omega \right \} \, .
\end{eqnarray*}

We equip $\Lambda_1$ with the norm
$$
\|f\|_{\Lambda_1} = \sup_{x, x+h \in \Omega \atop x \ne y}
  \frac{|f(x + h) + f(x - h) - 2f(x)|}{|h|}
   + \|f\|_{L^\infty(\Omega)} \, .
$$

Inductively, if $\alpha > 1$, we say
that $f \in \Lambda_\alpha(\Omega)$ if
$f$ is continuously differentiable, $f \in \Lambda_{\alpha-1}$,
and $\nabla f \in \Lambda_{\alpha - 1}$. The norm is
$$
\|f\|_{\Lambda_\alpha} = \|f\|_{\Lambda_{\alpha-1}(\Omega)} +
            \|\nabla f\|_{\Lambda_{\alpha-1}(\Omega)} \, .
$$

Assume that the domain $\Omega$ has $C^2$ boundary.
Let $U$ be a tubular neighborhood of
$\partial \Omega$ (see [HIR]).
Thus each point in $U$ has a unique nearest point in $\partial \Omega$
and there is a well-defined Euclidean orthogonal projection $\pi: U \ra \partial \Omega$.
We say that a $C^\infty$ curve $\gamma: [0,1] \ra \Omega \cap U$ lies
in ${\cal C}^1(\Omega)$ if {\bf (i)} $|\dot \gamma(t)| \leq 1$ for
all $t$ and {\bf (ii)} $\dot \gamma(t)$ lies in the complex
tangent space (see [KRA1]) at $\pi(\gamma(t))$ for each $t$. We think
of such a $\gamma$ as a ``normalized complex tangential curve''.
Let $0 < \alpha < \infty$.  Following E. M. Stein [STE2], we set
$$
\Gamma_{\alpha, 2\alpha}(\Omega) = \{f \in \Lambda_\alpha(\Omega):
 f \circ \gamma \in \Lambda_{2\alpha}([0,1]) \ \hbox{for each} \
 \gamma \in {\cal C}^1(\Omega)\} \, .
$$
Thus a function in $\Gamma_{\alpha, 2\alpha}(\Omega)$ is smooth
of order $\Lambda_\alpha$ in {\it all directions}, but smooth
of order $\Lambda_{2\alpha}$ in complex tangential directions.

In fact it is convenient to think about an $f \in \Gamma_{\alpha,2\alpha}$
as a function that is $\Lambda_\alpha$ along complex normal curves
and is $\Lambda_{2\alpha}$ along complex tangential curves. This
point of view has been developed, among other places, in
[KRA3]--[KRA6]. See also [GRS]  and [RUD].

Stein's remarkable theorem about these function spaces is as follows.

\begin{theorem} \sl
Let $\Omega \ss \CC^n$ be a domain with $C^2$ boundary.
Let $\alpha > 0$.
Let $f$ be a holomorphic function on $\Omega$ which lies
in $\Lambda_\alpha(\Omega)$. Then $f \in \Gamma_{\alpha, 2\alpha}(\Omega)$.
\end{theorem}

\noindent This result has been refined and
generalized in [KRA3]--[KRA6].

In [KRA2] it was shown that the analogous result for the space $BMO$
of functions of bounded mean oscillation fails.
To wit, if $\Omega =B \ss \CC^n$ is the unit ball,
then there are two types of balls to consider in the boundary $\partial B$:
If $P \in \partial B$ and $r > 0$ then
$$
\beta_1(P,r) = \{z \in \partial B: |z - P| < r\}
$$
and
$$
\beta_2(P,r) = \{z \in \partial B: |1 - z \cdot \overline{P}| < r\} \, .
$$
Notice that $\beta_1$ is the standard, isotropic Euclidean ball whereas
$\beta_2$ is a nonisotropic ball with extent $r$ in the complex normal
direction and extent $\sqrt{r}$ in the complex tangential directions.

We define
$$
BMO_1(\partial B) = \left \{ g \ \hbox{on} \ \partial B: \sup_{z, r} \frac{1}{\sigma(\beta_1(z,r))}
 \int_{\beta_1(z,r)} |g(\zeta) - g_{\beta_1(z,r)} | \, d\sigma(\zeta) < \infty \right \} \, .
$$
Here $d\sigma$ is  $(2n-1)$-dimensional Hausdorff measure
and $g_S$ denotes the average of $g$ over the set $S$: $g_S =
[1/V(S)] \cdot \int_S g(t) \, dV(t)$.
Likewise
$$
BMO_2(\partial B) = \left \{ g \ \hbox{on} \ \partial B: \sup_{z, r} \frac{1}{\sigma(\beta_2(z,r))}
 \int_{\beta_2(z,r)} |g(\zeta) - g_{\beta_2(z,r)} | \, d\sigma(\zeta) < \infty \right \} \, .
$$

Say that $g$ on the ball is in ${\cal BMO}_1(B)$ if $f$ is holomorphic, $f \in H^2$,
and $f$ has boundary function that is in fact in $BMO_1$.
Say that $g$ on the ball is in ${\cal BMO}_2(B)$ if $f$ is holomorphic, $f \in H^2$,
and $f$ has boundary function that is in fact in $BMO_2$.
The result of [KRA2] is that there is a
holomorphic function on the ball in $\CC^n$ which is in classical
isotropic ${\cal BMO}_1$ but not in nonisotropic ${\cal BMO}_2$. That proof used
quite a lot of functional analysis, and did not exhibit the
counterexample explicitly. A more concrete proof, with an explicit
example, was given in [ULR]. The result of [KRA2] was particularly
surprising because it showed as a byproduct that $BMO$ is not an interpolation
space between $L^p$ and $\Lambda_\alpha$.

In view of these results, it is natural to wonder whether a
holomorphic function on a domain $\Omega$ in $\CC^n$ that is
nontangentially bounded almost everywhere on a set $E \ss \partial \Omega$ will
in fact be admissibly bounded almost everywhere on $E$. If
this were true, then it would follow, at least for a
reasonable class of domains $\Omega$, that a function that is
nontangentially bounded on a set $E \ss \partial \Omega$ of
positive measure will in fact have admissible limits almost
everywhere on $E$. In view of the Lindel\"{o}f principle developed
in [CIK] and [KRA12], this is a very natural sort of result. And it turns
out to be true. Its proof is one of the main results of the present
paper.

\section{Definitions and Prior Results}

We take this opportunity to review the standard definitions and
concepts pertaining to this subject. The reference [KRA1] is
a good source for the details.  See also [KRA11].

We begin with harmonic analysis on $\RR^{N+1}$, which is the most natural
setting for the consideration of nontangential convergence.
Define the upper half space $$ U = \RR^{N+1}_+ = \{x = (x_1,
x_2, \dots, x_N, x_{N+1}):
  x_{N+1} > 0\} \, .
$$
We often shall write an element of $\RR^{N+1}_+$ as $(x', x_{N+1})$, where
$x' \in \RR^N$ and $x_{N+1} > 0$.  Of course the boundary of $\RR^{N+1}_+$
can be identified with $\RR^N = \{(x_1, \dots, x_N, 0)\}$ in a natural way.

If $\alpha > 1$ and $P \in \partial \RR^{N+1}_+ \equiv \RR^N$,
then we define the {\it Stolz region} or {\it nontangential
approach region} $\Gamma_\alpha$(P) of aperture $\alpha$ at $P$
to be
$$
\Gamma_\alpha(P) = \{x = (x', x_{N+1}) \in \RR^{N+1}_+: |x' - P| < \alpha x_{N+1}\} \, .
$$
This is a conical-shaped region in the upper half space. Points
in this region {\it cannot} approach the boundary along a
tangential curve.

Let $f$ be a function on $\rnp$. We say that $f$ is {\it nontangentially
bounded} on a set $E \ss \partial \rnp$ if, for each $P \in \partial
\rnp$, there is an $\alpha = \alpha(P) > 1$ such that
$f \bigr |_{\Gamma_\alpha(P)}$ is bounded. The bound, of course,
may (and, in general, will) depend on $P$ and on $\alpha$.
But observe that, if $E$ has positive $N$-dimensional measure, then we may use elementary
measure theory to find a set $E' \ss E$ of positive measure and
a constant $\alpha' > 1$ and a number $M' > 0$ so that
$|f| \leq M'$ on $\Gamma_{\alpha'}(P')$ for each $P' \in E'$.
Thus we may uniformize the estimate in the definition of ``nontangentially
bounded''.

With notation as in the last paragraph, we say that $f$ has
{\it nontangential limit} on the set $E \ss \partial \rnp$ if,
for each $\alpha > 1$, and each point $P \in E$, the limit
$$
\lim_{\Gamma_\alpha(P) \ni x \ra P} f(x)
$$
exists.\footnote{It is worth noting that the definition
of {\it nontangentially bounded} imposes on each $P \in E$ a
condition involving just one $\alpha$, depending on $P$.  But
the definition of {\it nontangential limit} imposes on each
$P \in E$ a condition for all $\alpha$.}

Calder\'{o}n's celebrated theorem [CAL] says this:

\begin{theorem} \sl
Let $u$ be a harmonic function on $\rnp$. Let $E \ss \partial \rnp$
have positive $N$-dimensional measure. If $u$ is nontangentially bounded
on $E$, then $u$ has nontangential limits almost everywhere on $E$.
\end{theorem}

Calder\'{o}n notes in his paper that his result holds for holomorphic
functions of several complex variables; but that more general
result is also formulated in terms of classical nontangential convergence.
See also [WID].  In fact the concept of admissible convergence would
not be invented for another twenty years.

Let $B \ss \CC^n$ be the unit ball. Let $f$ be a complex-valued
function on $B$ and let $E \ss \partial B$. We say that $f$ is
{\it admissibly bounded} on $E$ if, for each $P \in E$, there
is an $\alpha = \alpha(P) > 1$ such that $f$ is bounded on
${\cal A}_\alpha(P)$. Using elementary measure theory, it may be
seen (in analogy with the situation for classical nontangential convergence)
that if $E \ss \partial \Omega$ has positive
$(2n-1)$-dimensional measure then there is a set $E' \ss E$ of
positive measure and a number $\alpha' > 1$ and a constant
$M' > 0$ such that $|f|$ is bounded by $M'$ on ${\cal
A}_{\alpha'}(P')$ for each $P' \in E'$.  See [STE1].

With notation as in the last paragraph, we say that $f$ has
{\it admissible limit} on the set $E \ss \partial \Omega$ if,
for each $\alpha > 1$, and each point $P \in E$, the limit
$$
\lim_{{\cal A}_\alpha(P) \ni z \ra P} f(z)
$$
exists.

Now Stein's theorem [STE1, Theorem 12] states the following:

\begin{theorem} \sl
Let $\Omega \ss \CC^n$ be a strongly pseudoconvex domain
with $C^2$ boundary (see [KRA1]). Let $E \ss \partial \Omega$
be a set of positive $(2n-1)$-dimensional measure. Let
$f$ be a holomorphic function on $\Omega$. Then $f$ is
admissibly bounded on almost everywhere on $E$ if and only if $f$ has admissible
limits at almost every point of $E$.
\end{theorem}

The proof of this last result that appears in [STE1] relies on
the potential theory and the Levi geometry of the domain in
question. In particular, it requires the construction of a
special ``preferred'' Levi metric. It also depends on
estimates involving the Lusin area integral; that is to say,
the argument is not direct. S. Ross Barker [BAR] has provided
an alternative, more measure-theoretic approach to the matter
and thereby proved the result to be true on a broad class of
domains (and also avoided the use of the area integral).
Barker only enunciates and proves a result to the effect that
a holomorphic function that is admissibly bounded almost
everywhere (on the entire boundary) has then admissible limits
almost everywhere (on the entire boundary). He comments at the
end that his result can be localized (in the spirit of
Calder\'{o}n). One of the points of the present paper is to
provide a new approach to the Calder\'{o}n result. We can also prove a
sharper version of the theorem, in the sense that we can in
many cases adapt the shape of the approach regions to the Levi
geometry of the particular domain under study (see Section 8).

\section{The Main Results of the Present Paper}

In this section we collect the statements of the main results of the present
paper.  We also briefly indicate their context and significance.

Recall once again that, in the pioneering work [STE1], Stein proves
theorems about the boundary behavior of holomorphic functions {\it using
approach regions of the same parabolic complex tangential geometry}, no
matter what the particular intrinsic complex geometry of the domain in
question. It was only in later work (see [NSW1], [NSW2], [KRA9]) that the
mathematical machinery was developed for adapting the shape of the
approach region to the Levi geometry of the domain. The work in [DIB1],
[DIB2], [DIB3] extends the new ideas further. It should be stressed that the
results presented in the present paper build on these ideas. For
instance, the paper [BAR] certainly extends Stein's version of the
Calder\'{o}n local Fatou theorem to any smoothly bounded domain in
$\CC^n$; but it still used the old parabolic approach regions of Koranyi
and Stein. In the present paper we prove a version of this theorem for
several different types of domains; and we use approach regions that are
specifically adapted to the geometry of the domain in question (see
Sections 6, 8 for the details).

Our theorems do {\it not} apply to an arbitrary
smoothly bounded domain in $\CC^n$.
At this stage in the development of our mathematical machinery they
cannot. For all the proofs here require {\bf (i)} that the boundary of the
domain be equipped with a system of balls that, together with standard
$(2n-1)$-dimensional area measure, make the boundary a space of
homogeneous type in the sense of [COW1], [COW2]
and {\bf (ii)} the geometric
structure of the approach regions ${\cal A}_\alpha$ must be compatible (in
a sense to be described in detail below) with the balls from {\bf (i)}. As
of this writing, we know how to carry out such a program on {\bf (a)}
strongly pseudoconvex domains, {\bf (b)} domains of finite type in
$\CC^2$, and {\bf (c)} finite type, convex domains in $\CC^n$.  Refer
to Section 8 for the relevant geometric ideas.

\begin{theorem} \sl
Let $\Omega$ be either a strongly
pseudoconvex domain in $\CC^n$ or a finite type domain in
$\CC^2$ or a convex, finite type domain in $\CC^n$.
Let $E \ss \partial \Omega$ be a set of positive
measure (either 3-dimensional Hausdorff measure for a domain
in $\CC^2$ or $(2n-1)$-dimensional Hausdorff measure for a
domain in $\CC^n$). Suppose that $f$ is a holomorphic function
on $\Omega$. If $f$ is admissibly bounded at almost every
point of $E$ then $f$ has admissible limits at almost every
point of $E$.
\end{theorem}

\noindent {\bf Discussion:} In the strongly pseudoconvex case,
Stein proves this theorem in [STE1, Theorem 12]. His proof
proceeds by way of a Lusin area integral argument. We provide
a new, more direct proof and also extend the result to finite
type domains in $\CC^2$ and convex, finite type domains in
$\CC^n$. We avoid the use of a special ``preferred'' metric and of the
Lusin area integral and work more
directly with the Levi geometry of the domain.
As part of our treatment of Theorem 4, we shall need to
give a detailed consideration of approach regions for
Fatou theorems on domains of the type under discussion (Section 6).
This is a subtle matter, for the shape of the regions varies
in a sort of semi-continuous manner with the base boundary point.
Further details will also appear in Section 8 below.

\begin{theorem} \sl
Let $f$ be a holomorphic function on the unit ball in $\CC^n$, $n > 1$.
Let $E \ss \partial B$ be a set of positive $(2n-1)$-dimensional measure.
Then $f$ is nontangentially bounded at almost every point of $E$ if and
only if $f$ is admissibly bounded at almost every point of $E$.
\end{theorem}

\begin{corollary} \rm
Let $f$ be a holomorphic function on the unit ball in $\CC^n$, $n > 1$.
Let $E \ss \partial B$ be a set of positive $(2n-1)$-dimensional measure.
Assume that $f$ is nontangentially bounded at almost every point of $E$.
Then $f$ has admissible limits at almost every point of $E$.
\end{corollary}

\noindent {\bf Discussion:}  In fact this result is valid in considerably
greater generality.  But all the key ideas are already present in
the ball case, and matters are clearer when everything may be written
explicitly.

It should be stressed that Theorem 5 is {\it not} true point-by-point.
That is to say, at a particular point of the boundary of $B$ it is not
true that nontangential boundedness implies admissible boundedness.
This circle of questions is closely related to the Lindel\"{o}f principle,
for which see [CIK] and [KRA12].

The theorem answers a fairly old question, one that is rather natural
in view of the discussion in Section 1.  This new result puts the whole
idea of admissible convergence into a very natural context.

\section{An Ontology of Maximal Functions}

In this paper we shall use eleven different maximal functions.  For the
convenience of the reader, we collect all their definitions here.

We begin by thinking about the most natural and classical setting for
maximal functions, which is the Euclidean space $\RR^N$.  Let $f$ be
a locally integrable function on $\RR^N$.  For $x \in \RR^N$ we define
$$
Mf(x) = \sup_{r > 0} \frac{1}{|B(x,r)|} \int_{B(x,r)} |f(t)| \, dt \, .
$$
and
$$
{\cal M}f(x) = \limsup_{r \ra 0^+} \frac{1}{|B(x,r)|} \int_{B(x,r)} |f(t)| \, dt \, .
$$
Here, as usual,
\begin{itemize}
\item The set $B(x,r)$ is the standard isotropic Euclidean ball in $\RR^N$ with center $x$ and radius $r > 0$.
\item We let $|B(x,r)|$ denote the $N$-dimensional Lebesgue measure of $B(x,r)$, which
is $c_N r^N$.
\item The measure $dt$ is the standard Lebesgue measure.
\end{itemize}
It is also useful to let
$$
M_\delta f(P) = \sup_{0 < r \leq \delta} \frac{1}{|B(x,r)|} \int_{B(x,r)} |f(t)| \, dt \, .
$$

The first maximal operator $M$ is the classical one due to
Hardy and Littlewood. The second ${\cal M}$ is our first new maximal operator.
This maximal function differs from the classical one in that
the supremum has been replaced by the limit supremum.  The third maximal
operator $M_\delta$ is another small modification of $M$, restricting
to balls of radius not exceeding $\delta$.

Now let $\Omega \ss \CC^n$ be a domain on which a notion of admissible
approach region ${\cal A}_\alpha(P)$, $P \in \partial \Omega$, has been
defined---see Section 8.  If $g$ is a complex-valued
function on $\Omega$ and $P \in \partial \Omega$ then we define
$$
g_\alpha^{*}(P) = \sup_{z \in {\cal A}_\alpha(P)} |g(z)|
$$
and
$$
g_\alpha^{**}(P) = \limsup_{{\cal A}_\alpha(P) \ni z \ra P} |g(z)| \, .
$$

Now let $B \ss \CC^n$ be the unit ball.
Of course $\partial B$ is equipped with a family
of isotropic Euclidean balls
$$
\beta_1(P,r) = \{z \in \partial B: |z - P| < r\} \, .
$$
We shall also utilize the nonisotropic balls given by the condition
$$
\beta_2(P,r) = \{z \in \partial B: |1 - z \cdot \overline{P}| < r\} \, .
$$

Corresponding to these two types of balls in $\partial B$ we shall
have two types of maximal functions.  Let $d\sigma$ be boundary
area measure.  If $\varphi$ is a locally integrable
function on $\partial B$ and $P \in \partial B$ then we set
$$
M_1 \varphi(P) = \sup_{r > 0} \frac{1}{\sigma(\beta_1(P,r))} \int_{\beta_1(P,r)} |\varphi(\zeta)| \, d\sigma(\zeta)
$$
and
$$
M_2 \varphi(P) = \sup_{r > 0} \frac{1}{\sigma(\beta_2(P,r))} \int_{\beta_2(P,r)} |\varphi(\zeta)| \, d\sigma(\zeta)  \, .
$$
We also define two maximal functions based on the limsup rather than the supremum:
$$
{\cal M}_1 \varphi (P) = \limsup_{r \rightarrow 0^+} \frac{1}{|\beta_1(P,r)|} \int_{\beta_1(P,r)} |\varphi(t)| \, d\sigma(t)
$$
and
$$
{\cal M}_2 \varphi (P) = \limsup_{r \rightarrow 0^+} \frac{1}{|\beta_2(P,r)|} \int_{\beta_2(P,r)} |\varphi(t)| \, d\sigma(t) \, .
$$

In addition, we shall have two sets of truncated maximal operators as follows:
$$
M_{1,\delta} \varphi (P) = \sup_{0 < r \leq \delta} \frac{1}{\sigma(\beta_1(P,r))} \int_{\beta_1(P,r)} |\varphi(\zeta)| \, d\sigma(\zeta)
$$
and
$$
M_{2,\delta} \varphi (P) = \sup_{0 < r \leq \delta} \frac{1}{\sigma(\beta_2(P,r))} \int_{\beta_2(P,r)} |\varphi(\zeta)| \, d\sigma(\zeta) \, .
$$

\section{Some Estimates for Maximal Functions}

In this section we study some of our new maximal functions.  These functions
are rather natural tools for the study of the boundary behavior of
holomorphic functions.  Previous studies (see [STE1], [KRA1], [KRA9],
[NSW1], [NSW2]) endeavored to study the entire boundary at once---using
classical maximal functions that were designed for such a purpose.
Our goal here is to localize the process.  This change is particularly
propitious for the development of results of Calder\'{o}n type.

\begin{proposition} \sl
The maximal operator ${\cal M}$ is of weak type $(1,1)$ and also
of strong type $(p,p)$ for $1 < p \leq \infty$.
\end{proposition}
{\bf Proof:}  The classical Hardy-Littlewood maximal function $M$ is
known (see [STE3]) to be of weak type $(1,1)$ and of strong type $(p,p)$ for
$1 < p \leq \infty$.  Clearly ${\cal M} f(x) \leq Mf(x)$
for any $f$.  The result follows.
\endpf
\smallskip \\

We have formulated and proved Proposition 7 on $\RR^N$.  But the statement
and proof transfer {\it grosso modo} to the boundary of a $C^2$, bounded
domain in $\RR^N$ or $\CC^N$.  After all, such a boundary is a smooth
manifold hence is locally Euclidean.  Put in different terms, this boundary
is certainly a space of homogeneous type (see [COW1], [COW2])
when it is equipped with
isotropic balls and the standard Hausdorff measure on the boundary.  Observe that,
thus far, we are not taking the complex structure or the Levi geometry into account.
We are only looking at classical Euclidean geometry.

Our key tool in proving boundary limit theorems for holomorphic
functions is as follows.

\begin{theorem} \sl
Let $\Omega$ be a bounded domain in $\CC^n$, $n \geq 2$, with $C^2$ boundary which is of
one of these types:
\begin{itemize}
\item Strongly pseudoconvex domains in $\CC^n$;
\item Finite type domains in $\CC^2$;
\item Finite type, convex domains in $\CC^n$.
\end{itemize}
Let $u$ be a real-valued, nonnegative, plurisubharmonic function on $\Omega$, continuous on $\overline{\Omega}$.
Let $\varphi = u \bigr |_{\partial \Omega}$ be the boundary trace of $u$.
Let $P \in \partial \Omega$ and $\alpha > 1$.  Then
$$
u_\alpha^{**}(P) \leq C_\alpha \cdot {\cal M}_2 {\cal M}_1 \varphi (P) \, .  \eqno (\star)
$$
\end{theorem}

This is our local version of Lemma 8.6.10 in [KRA1] or Theorem
2, p.\ 11, the Lemma, p.\ 33, and Lemma 1b, p.\ 42 in [STE1].
It will be the key tool in obtaining a suitable version of
Calder\'{o}n's theorem for domains in $\CC^n$.
Note that the maximal functions on the righthand side of $(\star)$
are the new localized maximal functions defined in terms of
the limit supremum.  Thus the are {\it smaller} than the
maximal functions in the classical inequalities of Stein and Barker.

As we know,
once a result like Theorem 8 is established, then it is a
straightforward exercise with measure theory to see that
suitably bounded holomorphic functions have boundary limits.
What is new here is the local nature of the maximal function
estimate. This, coupled with the newly defined maximal
functions, will give a new way to think about Fatou-type
theorems and Calder\'{o}n-type theorems even in $\CC^1$.

We note that Theorem 8 has of course a standard, classical formulation
on the unit disc.  In that context, we deal of course with
nontangential convergence and there is only one
limsup-type maximal function
${\cal M}$ on the boundary (see [KRA1]).  The estimate then
reads
$$
u_\alpha^{**}(P) \leq C'' {\cal M} u(P) \, .
$$
Here we interpret $u_\alpha^{**}$ on the left to be the limsup
over $\Gamma_\alpha(P)$.

\section{Proof of Theorem 8 on the Disc and the Ball}

To fix ideas, we will begin by proving Theorem 8 on the disc $D$ in $\CC$.
Now of course the correct concept is classical nontangential convergence,
and the function $u$ is {\it subharmonic}.
Fix a point $P \in \partial D$.  We may as well suppose that $P = 1 + i0$.  Fix a
parameter $\alpha > 1$ and let $z \in \Gamma_\alpha(P)$ be near
the boundary.  Let $\delta = 1 -|z|$.
For a suitable $c > 0$, depending on $\alpha$, we may be sure that
$D(z, c\delta) \ss D$.  Thus
$$
u(z) \leq C \cdot M_\delta u (\pi(z)) \, ,
$$
where $\pi(z) = z/|z|$ is the standard Euclidean projection of $z$ to $\partial D$.
See [KRA1, Proposition 8.1.10] or [STE1] for the idea behind estimating $u$
by the classical maximal function on the boundary.

It is essential at this point to notice that the arc centered at $\pi(z)$ and
having radius $c'\delta$ will certainly contain $P$.  This is because
$z \in \Gamma_\alpha(P)$.  Note that $c' = c'(\alpha)$.  Hence we may estimate the last line by
$$
u(z) \leq C'' \cdot M_{c'\delta} u(P) \, .   \eqno (\star\star)
$$

Since we are now working on the unit disc $D$ in $\CC$, we no
longer have distinct maximal functions (modeled on the limsup)
based on either isotropic balls or nonisotropic balls.  There
is just the single limsup maximal function ${\cal M}$
based on arcs in $\partial D$.

Now choose a sequence $z_j \in \Gamma_\alpha(P)$ such that
$$
z_j \ra P \quad \hbox{and} \quad u(z_j) \ra \limsup_{\Gamma_\alpha(P) \ni z \ra P} u(z) \equiv u^{**}_\alpha(P) \, .
$$
Then we know by $(\star\star)$ that
$$
u(z_j) \leq C'' \cdot M_{c' \delta_j} u (P) \, ,
$$
Here $\delta_j = 1 - |z_j|$.  Certainly, since $z_j \ra P$, we know
that $\delta_j \ra 0$.  As $j \ra \infty$, the righthand side is
certainly $\leq C''' \cdot {\cal M} u(P)$.  We conclude therefore that
$$
u_\alpha^{**}(P) \leq C'' {\cal M} u(P) \, .
$$
That is the desired conclusion.
\endpf
\smallskip \\

Now let us turn to the situation on the ball $B \ss \CC^n$.
This circumstance is rather more delicate, for we cannot pass
directly from the interior to the boundary by way of single
maximal function in order to get the estimates that we need.  In
the end, the estimate that we obtain is in terms of {\it two}
maximal functions.

The nonisotropic balls mesh nicely with the admissible approach
regions
$$
{\cal A}_\alpha(P) = \{z \in B: |1 - z \cdot \overline{P}| < \alpha (1 - |z|)\} \, .
$$
Notice in particular that the set of points in ${\cal A}_\alpha(P)$ having
distance precisely $\delta > 0$ from $\partial B$ can be described by
$$
{\cal E}_\alpha(\delta) = \{z \in B: \delta_{\partial B}(z) = \delta,
    |1 - z \cdot \overline{P}| < \alpha \cdot \delta\} \, .
$$
The projection of ${\cal E}_\alpha(\delta)$ to $\partial B$ is
the set
$$
\{z \in \partial B: |1 - z \cdot \overline{P}| < \alpha \cdot \delta \} \, .
$$
Thus we see in a natural way that the admissible approach region
${\cal A}_\alpha(P)$ is built up from the nonisotropic balls
$\beta(P,r)$ and, conversely, the nonisotropic balls are projections
of level sets of the approach regions ${\cal A}_\alpha(P)$.  It is
this relationship, between balls and approach regions, that we shall
want to exploit when we study more general domains.

Another key ingredient of our analysis on the unit ball in
$\CC^n$ is the existence of certain polydiscs. If $\alpha > 1$
is fixed and $P \in \partial B$, then consider a point $z \in
{\cal A}_\alpha(P)$. It is helpful to normalize coordinates so
that $\Re z_1$ is the real normal direction at $z$ and thus
$\Im z_1$ is the complex normal direction. Thus $z_2, \dots,
z_n$ span the complex tangential directions at $z$. If
$\alpha' = 2\alpha$ then of course $z \in {\cal
A}_{\alpha'}(P)$. Thus, letting $\delta = 1 - |z|$, we see
that the polydisc
$$
{\cal D} = {\cal D}(z) \equiv D(z_1, \delta/2) \times D(z_2,, \sqrt{\delta/(2\alpha)}) \times \cdots \times
        D(z_2,, \sqrt{\delta/(2\alpha)})
$$
lies in ${\cal A}_{\alpha'}$ and hence in $B$.

Now, as usual, let $u$ be a nonnegative function that is continuous on $\overline{B}$ and
plurisubharmonic on $B$.  Certainly we have (iterating the sub-mean
value property on each coordinate disc in each dimension
that makes up ${\cal D}$)
$$
u(z) \leq \frac{1}{|{\cal D}|} \int_{\cal D} u(\zeta) \, dV(\zeta) \, .  \eqno (*)
$$
Now, in order to pass from the interior to the boundary, we must exploit
our knowledge of the classical Poisson integral.  Let us denote
the Poisson kernel by $P$ and the Poisson integral of a boundary
function $f$ by $Pf$.  It follows from the maximum principle that
the plurisubharmonic function $u$ is majorized by the Poisson
integral $P\varphi$
of its boundary function $\varphi$.  And that in turn is majorized by (see [KRA1, Chapter 8])
the Hardy-Littlewood maximal function $M_1 \varphi$ of $\varphi$ at the projected boundary
point of the argument; but we in fact only need the Hardy-Littlewood maximal
function based on balls of radius $\leq \delta$, and that we denote
by $M_{1, \delta}$.  Thus line $(*)$ is majorized by
$$
\frac{1}{|{\cal D}|} \int_{\pi(\cal D)} M_{1, \delta} \varphi(\pi(\zeta)) \cdot \delta \, d\sigma(\zeta) \, .
$$
Here, of course, $\pi(z) = z/|z|$ is the projection of
$B \setminus \{0\}$ to $\partial B$ and the extra $\delta$ in
the integrand comes from the real normal dimension of ${\cal D}$.

Now it is essential to note that $\pi({\cal D})$ is comparable to a nonisotropic
ball of center $\pi(z)$ and radius $c\delta$.  So we may
rewrite our estimate as
$$
u(z) \leq \frac{C}{\sigma(\beta(\pi(z), c\delta)} \int_{\beta(\pi(z), c\delta)} M_{1, \delta} \varphi(\zeta) \, d\sigma(\zeta) \, .
$$
Because the boundary is a space of homogeneous type---in particular the enveloping property is
valid---we may replace the ball of radius $c\delta$ and centered at $\pi(z)$ with a ball
of radius $c'\delta$ and centered at $P$.  So we have
$$
u(z) \leq \frac{C'}{\sigma(\beta(P, c'\delta))} \int_{\beta(P, c'\delta)} M_{1, \delta} \varphi(\zeta) \, d\sigma(\zeta) \, .
$$

And this line is not greater than
$$
C' M_{2, \delta}(M_{1, \delta} \varphi)(P) \, ,
$$
where $M_{2, \delta}$ is the Hardy-Littlewood-type maximal function modeled on the nonisotropic balls $\beta$ of
radius not exceeding $\delta$.
Now choose a sequence $z_j \in {\cal A}_\alpha(P)$ such that $u(z_j) \ra \limsup_{z \ra P} u(z)$.
Then of course
$$
u(z_j) \leq C' M_{2, \delta_j} (M_{1, \delta_j} \varphi)(P) \, .
$$
Letting $j \ra \infty$, we find that
$$
u^{**}_\alpha(P) \leq C' \cdot {\cal M}_2 {\cal M}_1 \varphi(P) \, .
$$
Of course the maximal functions ${\cal M}_j$ on the right denotes the ``limsup'' maximal
function that we defined and considered earlier.  Thus we have obtained the desired estimate.
\endpf
\smallskip  \\

We have only presented the proof of Theorem 8 so far on the unit ball $B$ in $\CC^n$.
But we assert that it is valid on more general classes of domains, as we have indicated
above.  In Section 8 we isolate those geometric properties that are needed in
order to see that the result goes through in the claimed greater generality.

\section{The (Localized) Calder\'{o}n Theorem}

Now we shall present our new approach to the Calder\'{o}n theorem.
To repeat, this point of view is new even in the classical setting
of the unit disc in $\CC$.  We shall confine our discussion
to the unit ball, where all the key ideas are already clear.

\begin{proposition}  \sl
Let $f$ be a holomorphic function on the unit ball $B \ss \CC^n$.
Let $M > 0$ and suppose that $|f| \leq M$.
Let $E \ss \partial B$ be a set of positive measure, and supposed
that $f$ is admissibly bounded on $E$.  Then $f$ has admissible
limits at almost every point of $E$.
\end{proposition}

\noindent {\bf Remark:}  The tauberian condition $|f| \leq M$ is a bit
artificial, and is certainly not part of the standard canon of
the Calder\'{o}n theorem.  But it is a useful tool in our proof.
Afterward, we shall remove this condition and recover the standard
Calder\'{o}n result.
\smallskip \\

\noindent{\bf Proof:} Let $\sigma$ be the usual rotationally-invariant
$(2n-1)$-dimensional area measure on $\partial B$. Let $\epsilon > 0$. By
outer regularity, select an open set $U \ss \partial B$ such that $U \supseteq E$ and
$\sigma(U \setminus E) \leq \epsilon \cdot \sigma(E)$. We shall use the
maximal functions, and the attendant notation, that we introduced earlier
in Section 4.

As usual, if $u$ is a plurisubharmonic function on $B$, continuous
on $\overline{B}$, and if $\varphi$ is the boundary trace of $u$, then
we know for each $P \in \partial B$ that
$$
u_\alpha^{**}(P) \leq C_\alpha {\cal M}_2 {\cal M}_1 \varphi (P) \, .
$$

Fix a point $P \in \partial \Omega$ and let $\nu = \nu_P$ be the unit
outward normal vector at $P$.
Following the classical argument presented in [KRA1, Theorem 8.6.11], we apply this
last inequality to the function
$$
f_{j,k}(z) = \left |f \left (z - \frac{1}{j}\nu_z\right ) - f\left (z - \frac{1}{k} \nu_z \right )\right |
$$
for some $j, k$ positive integers.
Then of course $f_{j,k}$ is plurisubharmonic.  If we restrict attention
to $f$ and $f_{j,k}$ on a neighborhood $\widetilde{\Omega} \cap \Omega$,
where $\widetilde{\Omega}$ is a neighborhood in $\CC^n$ of $P$, then we may also
take $f_{j,k}$ to be continuous on $\overline{\widetilde{\Omega} \cap \Omega}$.  Following
the argument in the proof of Theorem 8.6.10 in [KRA1], we know that
$$
\int_E |(f_{j,k})_\alpha^{**}(\zeta)|^2 \, d\sigma(\zeta) \leq
    C_\alpha \cdot \int_U {\cal M}_2 ({\cal M}_1 f_{j,k}(\zeta))^2 \, d\sigma(\zeta) \, .
$$
It is important to note that the maximal functions on the right are
defined using the limsup.  Thus we may say not only that each
maximal function is bounded on $L^2$, but also that it is bounded
from $L^2$ of any open set $\widetilde{U}$ containing $\overline{U}$
to $L^2$ of $U$---just because the boundedness would be proved using
arbitrarily small balls.  So we obtain
$$
\int_E |(f_{j,k})_\alpha^{**}(\zeta)|^2 \, d\sigma(\zeta) \leq
C''_\alpha \int_{\widetilde{U}} |f_{j,k}(\zeta)|^2 \, d\sigma(\zeta) \, ,
$$
where $\widetilde{U}$ is an open set in $\partial B$ that contains $\overline{U}$ and
such that $\sigma(\widetilde{U} \setminus E) < \epsilon \cdot \sigma(E)$.

Letting $j \ra \infty$ as in (8.6.10.2) of [KRA1], we find that
$$
\int_E \limsup_{{\cal A}_\alpha(\zeta) \ni z \ra \zeta} \left |(f\left (\zeta\right ) - f\left (z - \frac{1}{k} \nu\right )\right |^2 \, d\sigma(\zeta) \leq
C''_\alpha \int_{\widetilde{U}} \left |\widetilde{f}(\zeta) - f \left (\zeta - \frac{1}{k}\nu\right ) \right |^2 \, d\sigma(\zeta) \, .
$$
Here $\widetilde{f}(\zeta)$ denotes the nontangential limit of $f$ at almost every point of $E$, which
we know exists {\it a fortiori} by Calder\'{o}n's classical result.

Now of course the trick (on the righthand side) is to write $\widetilde{U} = (\widetilde{U} \setminus E) \cup E$.
Thus
$$
RHS = \int_{\widetilde{U} \setminus E} + \int_E  \equiv I + II \, .
$$
The first integral is estimated quite simply by $4M^2 \cdot \sigma(\widetilde{U} \setminus E) \leq C \cdot \epsilon \sigma(E)$.
Here, of course $C$ depends on $\alpha$ and on $M$.  But it does not depend on any of the other parameters that
are relevant to our present estimations.  In fact if we replace
$\epsilon$ by $\epsilon/M^2$, then we
may remove the dependence on $M$.  This will be important later.

So we have
\begin{eqnarray*}
\lefteqn{\int_E \limsup_{{\cal A}_\alpha(\zeta) \ni z \ra \zeta} \left |(f\left (\zeta\right ) - f\left (z - \frac{1}{k} \nu\right )\right |^2 \, d\sigma(\zeta)} \\
& \leq &
C''_\alpha \int_E \left |\widetilde{f}(\zeta) - f \left (\zeta - \frac{1}{k}\nu\right ) \right |^2 \, d\sigma(\zeta) + C \epsilon \sigma(E) \, . \qquad \qquad (*)
\end{eqnarray*}

And now one can proceed to imitate the argument at the end of the proof of Theorem 8.6.10 in [KRA1] to find that
$$
\sigma \left \{ \zeta \in E: \limsup_{{\cal A}_\alpha(\zeta) \ni z \ra \zeta} |f(z) - \widetilde{f}(\zeta)| > \epsilon \right \}
  \leq C \cdot \epsilon \cdot \sigma(E) \, .
$$
We conclude the argument with standard reasoning using elementary measure theory to see that
$\lim_{{\cal A}_\alpha(\zeta) \ni z \ra \zeta} f(z) = \widetilde{f}(\zeta)$.
\endpf
\smallskip \\

Our next job is to remove the tauberian hypothesis (i.e., the assumption of a global bound by $M$).

\begin{theorem}  \sl Let $f$ be a holomorphic function on the unit ball $B \ss \CC^n$.
Let $E \ss \partial B$ be a set of positive measure, and supposed
that $f$ is admissibly bounded almost everywhere on $E$.  Then $f$ has admissible
limits at almost every point of $E$.
\end{theorem}
{\bf Proof:}  For $\delta > 0$ small, let $B_\delta \equiv B(0, 1 - \delta) \subseteq B
\subseteq \CC^n$.  For each such $\delta > 0$ there is
of course a bound $M_\delta$ so that $|f| \leq M_\delta$ on $\overline{B_\delta}$.
If $E$ is as in the statement of the theorem, let $E_\delta$ be its Euclidean orthogonal projection
into $\partial B_\delta$.  Fix $\epsilon > 0$ as before.
Choose $\widetilde{U}_\delta \supset E_\delta$ so that $\sigma(\widetilde{U}_\delta \setminus E) < [\epsilon/M_\delta^2]\cdot \sigma(E)$.
Then the estimate $(*)$ holds on $B_\delta$ with
$E_\delta$ replacing $E$ (and, implicitly,
$\widetilde{U}_\delta$ replacing $\widetilde{U}$).
Now taking the supremum over $\delta > 0$,
we find that
\begin{eqnarray*}
\lefteqn{\int_E \limsup_{{\cal A}_\alpha(\zeta) \ni z \ra \zeta} \left |(f\left (\zeta\right ) - f\left (z - \frac{1}{k} \nu\right )\right |^2 \, d\sigma(\zeta)} \\
& \leq &
C''_\alpha \int_E \left |\widetilde{f}(\zeta) - f \left (\zeta - \frac{1}{k}\nu\right ) \right |^2 \, d\sigma(\zeta) + C \epsilon \sigma(E) \, .
\end{eqnarray*}
And now the proof may be completed as in the argument for the last theorem.
\endpf
\smallskip \\

A retrospective of the proof just presented shows that we have constructed machinery that
allows a standard sort of localization of the classical Fatou theorem.  If the ingredients
are in place to prove Theorem 8, then the Calder\'{o}n theorem follows immediately.
Section 8 explains how all these ingredients are present on domains other than the unit
ball $B$.

\section{Ingredients Needed for a Proof on a General Domain}

At this time we do not know how to prove the results considered here on a perfectly arbitrary
bounded domain in $\CC^n$ with $C^2$ boundary.  In fact our reasoning depends in essential ways (as does
the reasoning of Stein and others) on the Levi geometry of the domain.  The pertinent
desiderata are in fact known to hold on
\begin{enumerate}
\item[{\bf (i)}]  the unit disc in $\CC$;
\item[{\bf (ii)}]  the unit ball in $\CC^n$;
\item[{\bf (iii)}]  strongly pseudoconvex domains in $\CC^n$;
\item[{\bf (iv)}]  domains of finite type in $\CC^2$;
\item[{\bf (v)}]  convex domains of finite type in $\CC^n$, $n \geq 2$.
\end{enumerate}

We take this opportunity to isolate the essential features of the
geometry that are needed for our reasoning, and give references
where the reader may verify that these domains do indeed have the
required properties.  Fix a bounded domain $\Omega \ss \CC^n$ with
$C^2$ boundary.
\begin{enumerate}
\item[{\bf (a)}]  The boundary $\partial \Omega$ must be equipped with a family of balls
$\beta_2(P,r)$.  We use the notation $\beta_1(P,r)$ to denote the standard, isotropic,
Euclidean balls with center $P$ and radius $r$.  The ball $\beta_2(P,r)$ will typically
be nonisotropic and its shape will derive rather naturally from the complex structure
and/or the Levi geometry of $\Omega$.
\item[{\bf (b)}]  On the boundary of a suitable domain in $\CC^n$, the balls $\beta_2(P,r)$, together with the standard
$(2n-1)$-dimensional Hausdorff area measure $d\sigma$, form a space of homogeneous
type in the sense of [COW1], [COW2].   Of course the classical Euclidean balls
$\beta_1(P,r)$ together with $d\sigma$ also form a space of homogeneous
type.
\item[{\bf (c)}]  The domain $\Omega$ is equipped with a family of approach
regions ${\cal A}_\alpha(P)$ for each $P \in \partial \Omega$ and each $\alpha > 1$.
Each ${\cal A}_\alpha(P)$ is an open set in $\Omega$, and ${\cal A}_\alpha(P) \subseteq
{\cal A}_{\alpha'}(P)$ whenever $\alpha' > \alpha$.
\item[{\bf (d)}]  The approach regions ${\cal A}_\alpha(P)$ and the
balls $\beta_2(P,r)$ are related in the following manner.  If $\alpha > 1$ is
fixed and $\delta > 0$ is small then the Euclidean orthogonal projection
of
$$
\{z \in {\cal A}_\alpha(P): \delta_{\partial\Omega}(z) = \delta\}
$$
to $\partial \Omega$ is comparable to a ball $\beta_2(P, c\delta)$.  Here,
of course, $c$ will depend on $\alpha$.  Conversely, the set
$$
\bigcup_{\delta > 0} \{z \in \Omega: \pi(z) \in \beta_2(P, \delta),
   \delta_{\partial\Omega}(z) = \delta\}
$$
is comparable to an approach region ${\cal A}_{c\delta}(P)$.
\item[{\bf(e)}]  Suppose, after a normalization of coordinates, that
$\Re z_1$ is the real
normal direction at $z$, $\Im z_1$ the complex normal direction,
 and $z_2, \dots, z_n$ form an orthonormal basis
for the remaining $(n-1)$ complex tangential directions.
There is a $c > 0$ with the following property.
If $\alpha > 1$ is fixed and $z \in {\cal A}_\alpha(P)$
with $\delta = \delta_{\partial\Omega}(z)$ then
there are positive exponents $\lambda_1 = \lambda_1(z)$,
\dots, $\lambda_{n-1} = \lambda_{n-1}(z)$ so that the
polydisc
$$
{\cal D}(z) \equiv D(z_1, c\delta) \times D(z_2, c \delta^{\lambda_1}) \times \dots
 \times D(z_n, c\delta^{\lambda_{n-1}} )
$$
still lies in $\Omega$.
\item[{\bf (f)}]  A critical property of the polydisc ${\cal D}(z)$ in part {\bf (e)} is
that the Euclidean orthogonal projection $\pi({\cal D}(z))$ in $\partial \Omega$ is
comparable to a nonisotropic ball $\beta_2(\pi(z), c'\delta)$.  What is crucial
here is that $\delta$ will be the size of this ball in the complex normal direction,
and that will automatically determine all the other dimensions of the $(2n-1)$-dimensional
ball.
\item[{\bf (g)}]  The ball $\beta_2(\pi(z), c'\delta)$ from part {\bf (f)} is comparable
to a ball $\beta_2(P, c''\delta)$, where $P$ is as in part {\bf (e)}.
\end{enumerate}

A review of the proofs that we have presented in Sections 5, 6, 7 show that these seven
properties are precisely those that we used to establish our results.  Thus
Theorem
4 is true for the five types of domains described in {\bf (i)}--{\bf (v)}.

The references for properties {\bf (a)}--{\bf (g)} on domains {\bf (i)}--{\bf (iv)} are
\begin{enumerate}
\item[{\bf (i)}]  For the disc, see [KRA1].
\item[{\bf (ii)}]  For the ball, see [KRA11], [KRA1], [STE1].
\item[{\bf (iii)}]  For strongly pseudoconvex domains in $\CC^n$,
see [KRA1], [STE1], [KRL].
\item[{\bf (iv)}]  For finite type domains in $\CC^2$, see
[NSW1], [NSW2], [NRSW], [CAT].
\item[{\bf (v)}]  For convex, finite type domains in $\CC^n$, see
[DIF], [MCN1], [MCN2] and references therein.
\end{enumerate}

\section{The Nevanlinna Class}

For many purposes, the most natural space of functions on which to
consider Fatou-type theorems is the Nevanlinna class.  Here,
for a fixed bounded domain $\Omega \ss \CC^n$ with $C^2$ boundary, we
say that $f$ on $\Omega$ lies in ${\cal N}^+$ if {\bf (i)} $f$ is
holomorphic and {\bf (ii)} $\log^+ |f|$ has a harmonic majorant.
By a standard lemma that can be found in [STE1] or [KRA1], this
definition is equivalent to requiring that
$$
\sup_{0 < \epsilon < \epsilon_0} \int_{\partial \Omega_\epsilon}
\log^+ |f(\zeta)| \, d\sigma(\zeta) < \infty \, .
$$
Here $\Omega_\epsilon = \{z \in \Omega: \rho(z) = - \epsilon\}$
for some defining function $\rho$ for $\Omega$ (see [KRA1]) and
$$
\log^+ x = \left \{ \begin{array}{lcr}
               0 & \hbox{if} & x \leq 1 \\
             \log x & \hbox{if} & x > 1 \, .
            \end{array}
       \right.
$$

Stein's book [STE1] contained rather elaborate and technical arguments
to handle the boundary behavior of functions in ${\cal N}^+$.  A few
years later, Barker [BAR] provided a much simpler approach.  His key ideas
was the next lemma.  Note also that the case of {\it meromorphic functions}
in the Nevanlinna class was treated by Neff [NEF1], [NEF2] and Lempert [LEM].

\begin{lemma} \sl
Let $u$ be a nonnegative, continuous, plurisubharmonic function on $\Omega$
(we do not necessarily mandate that $u$ be continuous on $\overline{\Omega}$).
 Assume that $u$ has a harmonic majorant.  [Thus there is a finite, positive
 measure $\mu$ on $\partial \Omega$ such that
$$
u(z) \leq \int_{\partial \Omega} P(z, \zeta) \, d\mu(\zeta) \, .]
$$
Here of course $z \in \Omega$ and $P$ is the standard Poisson kernel.  Let $\alpha > 1$.
Then the admissible maximal function
$$
u^{*}_\alpha (\zeta) \equiv \sup_{z \in {\cal A}_\alpha(\zeta)} |u(z)|
$$
for $\zeta \in \partial \Omega$ satisfies
$$
u^{*}_\alpha (\zeta) \leq C_\alpha \left [ M_2( [M_1(\mu)]^{1/2}) \right ]^2 \, .
$$
and hence is finite almost everywhere in $\partial \Omega$.
\end{lemma}

We note first of all that Barker's lemma is still true if we replace $u^*_\alpha$ with
our maximal function $u^{**}_\alpha$ (defined using the
limsup), $M_1$ with ${\cal M}_1$, and
$M_2$ with ${\cal M}_2$.  Thus we know that
$$
u^{**}_\alpha (\zeta) \leq C_\alpha \left [ {\cal M}_2( [{\cal M}_1(\mu)]^{1/2}) \right ]^2 \, .
$$

As Barker notes, in case $f \in {\cal N}^+$,
one may apply this last lemma to the function $u = \log^+|f|$.  It follows
then that $u^{**}_\alpha$ is finite almost everywhere, and we may then
use our standard arguments to see that $f$ has an admissible limit almost everywhere.
Thus we have

\begin{theorem} \sl
Let $\Omega \ss \CC^n$, $n \geq 2$, be a bounded domain
with $C^2$ boundary. Assume that either $\Omega$ is the ball,
or a finite type domain in $\CC^2$, or a convex finite type
domain in $\CC^n$. Suppose that $f \in {\cal N}^+(\Omega)$.
Then $f$ has admissible boundary limits almost everywhere.
\end{theorem}

\section{Nontangential Versus Admissible Approach}

Now we shall prove Theorem 5. In fact, following the example that we have
already set with our proof of Theorem 4 (see Proposition 9), we shall at first
prove a version of the theorem that has an additional tauberian
hypothesis.

\begin{theorem} \sl
Let $f$ be a holomorphic function on the
unit ball $B$ in $\CC^n$, $n > 1$. Assume that there is a constant
$M > 0$ so that $|f| \leq M$. Let $E \ss \partial B$ be a set
of positive $(2n-1)$-dimensional measure. Then $f$ is
nontangentially bounded at almost every point of $E$ (with a
bound $C$ that is in general, and most interestingly, smaller
than $M$) if and only if $f$ is admissibly bounded (with the
same bound $C$) at almost every point of $E$.
\end{theorem}

As enunciated, we shall work on the domain the ball $B$, and for simplicity
and clarity we shall restrict attention to $B \ss \CC^2$. Thus assume that
the holomorphic function $f$ on $B$ is nontangentially bounded on the set
$E \ss \partial B$ of positive 3-dimensional Hausdorff measure. As usual we
call the measure $d\sigma$.

With elementary measure-theoretic arguments, we may extract from
$E$ a subset of positive measure so that $f$ is nontangentially bounded
at each point of the subset with a {\it uniform bound} $C$ and on
a cone $\Gamma_\alpha$ of uniform size---independent of the point.
We continue to call this new set $E$.  Not that, in general,
$C < M$---that is certainly the most interesting case.
We shall show then that $f$ is {\it admissibly bounded}
with bound $C$.

Now let $P \in \partial B$ be a point of density
(with respect to classical, isotropic balls)
of $E \ss \partial B$, and
let $U \ss \partial B$ be a small, relatively open neighborhood of $P$.
Let us consider a foliation of $U$ by complex tangential curves.
Call the curves $\gamma_w: (-\epsilon, \epsilon) \rightarrow U$, where $w$
is a 2-dimensional parameter.  Let $g_w$ denote the image curve of $\gamma_w$.
Restrict attention now to those
$g_w$ which intersect $E$ in a set of positive 1-dimensional measure.
For each such $g_w$, pick a point $\gamma_w(t_w)$ that is
a point of 1-dimensional density of $g_w \cap E$.   Let $\epsilon > 0$.
Choose a neighborhood $I_w = (t_w - \delta_w, t_w + \eta_w)$ so that
$$
\frac{{\cal H}^1(\gamma_w(I_w) \cap E)}{{\cal H}^1(I_w)} > 1 - \epsilon \, .
$$
We may suppose that $t_w, \delta_w, \eta_w$ are rational numbers.  Now, with some
elementary measure theory, we may focus on a collection of $\gamma_w$,
$w$ in a 2-dimensional set of positive measure, so that each of the
$I_w$ is the {\it same} interval $I^*$.  Give this set of $w$ the name ${\cal S}$,
and let $s \in {\cal S}$ be a 2-dimensional point of density.  We fix attention
on the point $x_0 = \gamma_s(t_s)$.

We may repeat the preceding arguments using a foliation $\widetilde{\gamma}_w$ of $U$
that is still complex tangential but is {\it transverse} to $\gamma_w$ (remember
that we are working in the boundary of the ball $B$ in $\CC^2$, so the complex
tangent space has real dimension 2).  This gives rise to a point $\widetilde{s} \in \widetilde{S}$.
By elementary measure theory---in particular by Fubini's theorem---we may suppose
that $x_0 = \gamma_s(t_s) = \widetilde{\gamma}_{\tilde{s}}(\widetilde{t}_{\tilde{s}}) = \widetilde{x}_0$.   We continue to call the point $x_0$.

Thus we focus our attention on the curves $\gamma_w(I^*)$ for $w
\in {\cal S}$ and $\widetilde{\gamma}_{\widetilde{w}}(\widetilde{I}^*)$ for $\widetilde{w} \in \widetilde{S}$.
We examine an admissible approach region with base point $x_0$ as above.
Call that region ${\cal A}_{\alpha}(x_0)$, some $\alpha > 1$.
Let $z \in {\cal A}_\alpha(x_0)$ be near
to the boundary---at distance much less than the length of
$\widetilde{I}$ or $\widetilde{I}^*$. Let $\delta = \delta_{\partial B}(z)$. Now consider, as
usual, a nonisotropic polydisc ${\cal D}$ centered at $z$, having radius $c'\delta$
in the complex normal directions and radii $c'\sqrt{\delta}$ in the
complex tangential directions, some small $c' > 0$.

The natural thing to do at this point is to estimate
$$
|f(z)| \leq \frac{1}{|{\cal D}|} \int_{{\cal D}} |f(\zeta)| \, dV(\zeta) \, .
$$
Because of our density statements about $\widetilde{I}$ and ${\cal S}$, we
can estimate this last line by
$$
  (1 - c''\epsilon) C + c'' \epsilon \cdot M \, .
$$
Since the point $z \in {\cal A}_\alpha(s)$ was chosen arbitrarily,
and since $\epsilon > 0$ was arbitrary, we in fact have shown
that $f$ is admissibly bounded at $x_0$ with bound $C$.  Since points
of the kind $x_0$ are measure-theoretically generic, we now know that we have a set of
positive measure in $E$ on which $f$ is admissibly bounded.
Again, by elementary measure theory, we may then conclude
that $f$ is admissibly bounded at almost all points of $E$.
That completes the proof.
\endpf
\smallskip \\

It remains to show that our result holds without the tauberian
hypothesis $|f| \leq M$.  So now let $f$ be nontangentially bounded
on a set $E \ss \partial B$ of positive measure.
As usual, we may take the nontangential approach regions $\Gamma_\alpha(P)$
to be of uniform aperture, and the bound $C$ to be uniform.

For $\tau > 0$ small, let $B_\tau = B(0,1 - \tau)$. Then of course $f$ is bounded
by some $M_\tau$ on $B_\tau$. Let $E_\tau$ be the projection of $E$ to
$\partial B_\tau$.  Then of course $f$ is nontangentially bounded on $E_\tau$
by $C$ (because each approach region
${\cal A}_\alpha^\tau(P_\tau) \subseteq B_\tau$ for
$P_\tau \in E_\tau$ is a subset of ${\cal A}_\alpha(P)$, where
$P = \pi(P_\tau)$).  Since the tauberian hypothesis is in place on $B_\tau$, we may
conclude that $f$ is admissibly bounded by $C$ on $E_\tau$.  But now,
for each $P \in E$, note that
$$
{\cal A}_\alpha (P) = \bigcup_{\tau > 0 \ \rm small} {\cal A}_\alpha^\tau (P_\tau) \,
$$
where ${\cal A}_\alpha^\tau(T_\tau)$ is the admissible region in $B_\tau$
based at the point $P_\tau$ (the projection of $P$ to $\partial B_\tau$).
Since $f$ is admissibly bounded by $C$ on each of the approach regions
on the right, it follows that $f$ is bounded by $C$ on ${\cal A}_\alpha(P)$.
This reasoning is valid at almost every point $P$ of $E$.  The proof
is therefore complete.

\section{Concluding Remarks}

The results in this paper are formulated and proved on the ball,
on strongly pseudoconvex domains, on finite type domains in $\CC^2$,
and on convex, finite type domains in $\CC^n$.  Other types
of domains can be handled with {\it ad hoc arguments}.  Among those
are the bidisc and complete Reinhardt domains like
$$
\Omega_{2, \infty} = \{z \in \CC^2: |z_1|^2 + 2 e^{-1/|z_2|^2} < 1\} \, .
$$
A complete theory of Fatou theorems and Calder\'{o}n theorems, which
can treat any bounded $C^2$ domain and which fully accounts for its
attendant Levi geometry, has yet to be produced.  The paper [KRA9] offers
a conceptual framework for handling all domains---using the Kobayashi metric
as a stepping stone and structural tool---but in practice it is rather
difficult to verify all the hypotheses of the results in [KRA9].

We are of the opinion, however, that invariant metrics are the right argot
for formulating function theoretic problems and results on arbitrary
domains.  Such metrics can read the Levi geometry, and they also take
into account the way that holomorphic functions in the interior depend
on the shape of the domain.  We look forward to future work in this
direction.

\newpage

\noindent {\LARGE \sc References}
\vspace*{.2in}

\begin{enumerate}

\item[{\bf [BAR]}]  S. R. Barker, Two theorems on boundary values
of analytic functions, {Proc.\ Am.\ Math.\ Soc.} 68(1978),
54--58.

\item[{\bf [CAL]}]  A. P. Calder\'{o}n, On the behavior of harmonic
functions near the boundary, {\em Trans.\ Am.\ Math.\ Soc.} 68(1950),
47-54.

\item[{\bf [CAT]}]  D. Catlin, Estimates of invariant metrics on pseudoconvex
domains of dimension two, {\it Math.\ Zeit.} 200(1989), 429--466.

\item[{\bf [CIK]}]  J. A. Cima and S. G. Krantz, A Lindel\"{o}f principle and
normal functions in several complex variables, {\em Duke Math.\ Jour.}
50(1983), 303-328.

\item[{\bf [COW1]}] R. R. Coifman and G. Weiss, {\it Analyse Harmonique
Non-Commutative sur Certain Espaces Homogenes}, Lecture Notes in Math.\
242, Springer-Verlag, Berlin, 1971.

\item[{\bf [COW2]}]  R.\ R.\  Coifman and G.\ Weiss, Extensions of
Hardy spaces and their use in analysis, {\it Bull.\ AMS}
83(1977), 569-645.

\item[{\bf [DIB1]}] F. Di Biase, Approach regions and maximal functions in
theorems of Fatou type, thesis, Washington University in St.~Louis, 1995.

\item[{\bf [DIB2]}] F. Di Biase, Exotic Convergence in Theorems of Fatou
Type, in {\it Harmonic Functions on Trees and Buildings}, Adam Koranyi,
ed., Contemporary Mathematics, vol.\ 206, American Mathematical Society,
1997.

\item[{\bf [DIB3]}]  F. Di Biase, {\it Fatou type theorems:
Maximal Functions and Approach Regions}, Birkh\"{a}user Publishing,
Boston, 1998.

\item[{\bf [DIF]}]  F. Di Biase and B. Fischer,
Boundary behaviour of $H^p$ functions on convex domains
of finite type in $\CC^n$, {\it Pacific J.\ Math.} 183(1998), 25--38.

\item[{\bf [DIK]}]  F. Di Biase and S. G. Krantz, {\it The Boundary Behavior
of Holomorphic Functions}, Birkh\"{a}user Publishing, Boston, 2006, to appear.

\item[{\bf [FAT]}]  P. Fatou, S\'eries trigonom\'etriques et s\'eries
de Taylor, {\it Acta Math.}, 30(1906), 335-400.

\item[{\bf [GAR]}]  J. B. Garnett, {\it Bounded Holomorphic Functions},
Academic Press, New York, 1981.

\item[{\bf [GRS]}] P. Greiner and E. M. Stein, {\it Estimates for the
$\overline{\partial}$-Problem}, Princeton University Press, Princeton,
NJ, 1977.

\item[{\bf [HIR]}]  M. Hirsch, {\it Differential Topology}, Springer-Verlag,
New York, 1976.

\item[{\bf [HUW]}]  R. Hunt and R. Wheeden,  On the boundary
values of harmonic functions 132(1968), 307--322.

\item[{\bf [JEK]}]  D. Jerison and C. Kenig, Boundary behaviour of
harmonic functions in non-tangentially accessible
domains, {\it Adv.\ Math.} 46(1982), 80--147.

\item[{\bf [KOR1]}] A. Koranyi, Harmonic functions on Hermitian hyperbolic
space, {\em Trans.\ A. M. S.} 135(1969), 507-516.

\item[{\bf [KOR2]}]  A. Koranyi, Boundary behavior of Poisson integrals
on symmetric spaces, {\em Trans.\ A.M.S.} 140(1969), 393-409.

\item[{\bf [KRA1]}]  S. G. Krantz, {\it Function Theory of Several Complex Variables},
$2^{\rm nd}$ ed., American Mathematical Society, Providence, RI, 2001.

\item[{\bf [KRA2]}] S. G. Krantz, Holomorphic functions of
bounded mean oscillation and mapping properties of the
Szeg\"{o} projection, {\em Duke Math. J.} 47(1980), 743-761.

\item[{\bf [KRA3]}] S. G. Krantz, Intrinsic Lipschitz classes on manifolds with applications
to complex function theory and estimates for the  $\dbar$
and  $\dbar_{b}$  equations,  {\it Manuscripta Math.} 24(1978), 351-378.

\item[{\bf [KRA4]}]  S. G. Krantz, Smoothness of harmonic and holomorphic functions,
{\em Proc.\ Symp.\ Pure Math.}, Vol. 35 (1979)
(S. Wainger and G. Weiss, eds.),    63-67.

\item[{\bf [KRA5]}]  S. G. Krantz, Characterizations of various domains of holomorphy
via $\dbar-$ estimates and applications to a problem of Kohn, {\it
Illinois J. Math.} 23(1979), 267-285.

\item[{\bf [KRA6]}]  S. G. Krantz, Lipschitz spaces on stratified groups,  {\it Trans.\  Am.\
Math. Soc.} 269(1982), 39-66.

\item[{\bf [KRA7]}]  S. G. Krantz, Finite type conditions and elliptic boundary
value problems, {\it Jour.\ Diff.\ Eq.} 34(1979), 239-260.

\item[{\bf [KRA8]}]  S. G. Krantz, Estimation of the Poisson kernel, {\it Journal
of Math.\ Analysis and Applications} 302(2005), 143--148.

\item[{\bf [KRA9]}]  S. G. Krantz, Invariant metrics and the boundary behavior of holomorphic
functions on domains in $\CC ^{n},$ {\em Jour.\ Geometric Anal.} 1(1991), 71-98.

\item[{\bf [KRA10]}]  S. G. Krantz, Fatou theorems old and new:  an overview of the
boundary behavior of holomorphic functions, Proceedings of
an International Conference on Complex Variables held in Seoul, Korea,
{\it Journal of the Korean Math.\ Society}37(2000), 139--175.

\item[{\bf [KRA11]}]  S. G. Krantz, {\it Partial Differential Equations and
Complex Analysis}, CRC Press, Boca Raton, FL, 1992.

\item[{\bf [KRA12}]  S. G. Krantz, The Lindel\"{o}f principle in
several complex variables, preprint.

\item[{\bf [KRL]}]  S. G. Krantz and S.-Y. Li, Area integral
characterizations of Hardy spaces on domains in $\CC^n$,
{\it Complex Variables} 32(1997), 373--399.

\item[{\bf [LEM]}] L. Lempert, Boundary behavior of meromorphic
functions of several complex variables, {\em Acta Math.}
144(1980), 1-26.

\item[{\bf [MCN1]}]  J. McNeal, Convex domains of finite type, {\it J. Funct.\
Anal.} 108(1992), 361--373,

\item[{\bf [MCN2]}]  J. McNeal, Estimates on the Bergman kernels
of convex domains, {\it Advances in Math.} 109(1994),
108--139.

\item[{\bf [NRSW]}] A. Nagel, J.-P. Rosay, E. M. Stein, and S. Wainger,
Estimates for the Bergman and Szeg\"o kernels in $\CC^2$, {\it Ann.\ of
Math.} 129(1989), 113--149.

\item[{\bf [NSW1]}]  A. Nagel, E. M. Stein, and S. Wainger, Boundary
behavior of functions holomorphic in domains of finite type, {\em Proc.\
Nat.\ Acad.\ Sci.\ USA} 78(1981), 6596-6599.

\item[{\bf [NSW2]}]  A. Nagel, E. M. Stein, and S. Wainger, Balls and
metrics defined by vector fields I: Basic properties, {\it Acta Math.}
155(1985), 103-147.

\item[{\bf [NEF1]}]  C. A. Neff, Maximal Function Estimates for Meromorphic Nevanlinna Functions, thesis,
Princeton University, 1986.

\item[{\bf [NEF2]}]  C. A. Neff, Boundary Convergence of Functions in the Nevanlinna Class, {\it Colloq.\ Math.}
60(1990), 477-506.

\item[{\bf [PLE]}]  A. Plessner, \"{U}ber die Verhalten analytischer
Funktionen am Rande ihres Definitionsbereiches, {\it J. F. M.}
159(1927), 219--227.

\item[{\bf [PRI1]}]  I. I. Privalov, Integrale de Cauchy,
Saratov, 1919.

\item[{\bf [PRI2]}] I. I. Privalov, {\it Randeigenschaften
analytischer funktionen}, $2^{\rm nd}$ ed., VEB Deutscher
Verlag der Wissenschaften, Berlin, 1956.

\item[{\bf [RUD]}] W. Rudin, Holomorphic Lipschitz functions in
balls, {\em Comment.\ Math.\ Helvet.} 53(1978), 143-147.

\item[{\bf [STE1]}] E. M. Stein, {\it Boundary Behavior of
Holomorphic Functions of Several Complex Variables}, Princeton
University Press, Princeton, 1972.

\item[{\bf [STE2]}] E. M. Stein, Singular integrals and
estimates for the Cauchy-Riemann equations, {\em Bull.\ A.M.S.}
79(1973), 440-445.

\item[{\bf [STE3]}]  E. M. Stein, {\it Singular Integrals and
Differentiability Properties of Functions}, Princeton University Press,
Princeton, NJ, 1970.

\item[{\bf [TSU]}]  M. Tsuji, bdry beh. of holo. fcns., 1930s  xxxxxx

\item[{\bf [ULR]}]  D. Ullrich, Tauberian theorems for pluriharmonic functions
which are BMO or Bloch, {\it Mich.\ Math.\ Jour.} 33(1986), 325-333.

\item[{\bf [WID]}]  K. O. Widman, On the boundary behavior of solutions to
a class of elliptic partial differential equations, {\it Ark.\ Mat.}
6(1966), 485--533.

\item[{\bf [ZYG1]}]  A. Zygmund, On the boundary values of functions
of several variables. I, {\it Fund.\ Math.} 36(1949), 207--235.

\item[{\bf [ZYG2]}]  A. Zygmund, A remark on functions of several variables,
{\it Acta Sci.\ Math.\ Szeged}, Pars B, 12(1950), 66--68.

\item[{\bf [ZYG3]}]  A. Zygmund, Note on the boundary values of functions
of several variables, {\it Ann.\ of Math. Studies} 25, Princeton
University Press, Princeton, NJ, 1950.

\end{enumerate}

\end{document}